\documentclass{article}

\def\R{{\rm I\! R}}

\newtheorem{theorem}{Theorem}

\newtheorem{lemma}{Lemma}

\newtheorem{remark}{Remark}

\title{Normalizers of planar systems with known first integrals}

\author{M. Sabatini}
\date{}
\begin{document}
\maketitle
\begin{abstract}  Given a planar differential system with a first integral, we show how to find a normalizer. For systems with a center, we give an integral formula for the derivative of its period function. 
\footnote[0]{
Partially supported by the PRIN group \lq\lq Equazioni differenziali ordinarie e applicazioni \rq\rq and by the GNAMPA group  \lq\lq Analisi qualitativa per sistemi dinamici in dimensione finita \rq\rq.}
\end{abstract}

\section{Introduction}

Let us consider a planar differential system
\begin{equation}\label{sysV} 
z' = V(z), \qquad z\in \Omega \subset \R^2,
\end{equation} 
with $\Omega$ open connected, $V(z)=  (V_1(z),V_2(z))\in C^2(\Omega,\R^2)$, $z = (x,y) \in \Omega$. 

A connected subset $A$ of $\Omega$ is said to be a {\it period annulus} of (\ref{sysV}) if every orbit of $V$ contained in $A$ is a non-trivial cycle of (\ref{sysV}). If the inner boundary of $A$ consists of a single point $O$, then  $O$ called a {\it center}, and the largest connected punctured neighbourhood $N_O$ of $O$ covered with non-trivial cycles is its {\it central region}. If $A$ is a period annulus, we can define on $A$ the {\it period function T} by assigning to each point $z\in A$ the minimal period $T(z)$ of the cycle $\gamma_z$ passing through $z$. We say that the period function $T$ is {\it increasing} if external cycles have larger periods. $A$ is said to be {\it isochronous} if $T$ is constant on $A$. 
Let us consider a curve $\delta(s)$ of class $C^1$ meeting transversally the cycle $\gamma$ at the point $s=s_0$. We say that $\gamma$ is a {\it critical cycle} if $\left[ \frac{d}{ds}T(\delta(s))\right]_{s=s_0} =0 $. It is possible to prove that such a definition does not depend on the particular transversal curve $\delta$ chosen. 
 
One is lead to study $T$'s monotonicity while approaching several problems related to (\ref{sysV}), as  boundary value problems, bifurcation or perturbation problems (\cite{Ch2}, \cite{Sch}). Moreover, it appears also also in relation to delay differential equations  \cite{CGM}, thermodynamics (\cite{R1}, \cite{R2}), linearizability \cite{MRT}. Finally, isochronicity is strictly related to stability, since a periodic solution contained in $A$ is Liapunov stable if and only if it has an isochronous neighbourhood \cite{CS}.

Recently,  $T$'s monotonicity has been studied by means of a sepcial class of auxiliary planar systems, called {\it normalizers}.   
Given a second vector field $U(z)=  (U_1(z),U_2(z))\in C^2(\Omega,\R^2)$, let us denote by $[V,U] = \partial_V U -\partial_U V$ the Lie bracket of $V$ and $U$. If $U$ is transversal to $V$, then it 
is said to be a {\it non-trivial normalizer} of $V$ on a set $A\subset \Omega$ if $[V,U] \wedge V =  0$ on $A$. If $U$ is a normalizer of $V$ on $A$, there exists a $C^1$ function $\mu$, defined on $A$, such that $[V,U] = \mu V$, where $\mu =\frac{ [V,U] \cdot  V}{|V|^2}$. The algebraic property $[V,U] \wedge V =  0$ has a simple dynamical interpretation.  Let $\phi_V(t,z)$ and $\phi_U(s,z)$ be the local flows defined by the solutions of (\ref{sysV}), and, respectively, by the solutions of
\begin{equation}\label{sysU} 
z' = U(z) .
\end{equation} 
In \cite{FGG}, theorem 1, it was proved that the derivative of $T$ along the solutions of (\ref{sysU}) is given by the following formula,
\begin{equation}\label{formulaFGG}
\partial_U T(z) = \left[ \frac{d}{ds} T(\phi_U(s,z)) \right]_{s=0}= \int_0^T \mu(\phi_V(t,z)) dt =
 \int_0^T \mu(\gamma_z(t)) dt.
\end{equation} 

In the same paper a non-trivial normalizer was found for Hamiltonian systems with separable variables and the monotonicity of the period function was studied in detail for centers of such systems. 

Such an approach was limited by the necessity to find a normalizer. This motivated the results of  \cite{S1}, where a normalizer was replaced by a transversal vector field $W$, proving the following formula for $T$'s derivative   
\begin{equation} \label{derivata}
\partial_W T(z) = \frac{1}{\beta(z)} \int_0^T \eta\big(\gamma_z(\tau)\big) \beta\big(\gamma_z(\tau)\big) \ d\tau .
\end{equation} 
In the above formula $\eta$ is a known function, 
\begin{equation}\label{etanu}
\eta = \frac{[V,W] \wedge W}{V \wedge W},
 \end{equation} 
 while $\beta$ is only implicitly known, since it is defined by means of an integral involving the known function $\nu$,
 \begin{equation} \label{intbeta}
 \nu = \frac{[V,W] \wedge V}{W \wedge V}, \qquad
\beta(\gamma_z(t)) = \beta(\gamma_z(0))\exp \left( -\int_0^t \nu(\gamma_z(\tau))\ d\tau \right).
\end{equation}
Even if this alternative approach allows to avoid some of the constraints related to the method introduced in \cite{FGG}, its application is bounded by the fact that $\beta$ is unknown.
 
In this paper we turn back to the approach of  \cite{FGG}, constructing explicitly a normalizer for a planar system with a first integral $H$. 
We first do it for hamiltonian systems and compute the corresponding function $\mu$, then we  pass to non-hamiltonian systems giving a normalizer and its $\mu$. The normalizer we find,
\begin{equation}\label{normH}
x' = \frac{H_x}{|\nabla H|^2} ,  \qquad y' = \frac{H_y}{|\nabla H|^2},
\end{equation}
is such that the derivative of $T$ along its solutions is just $T'(H)$. Moreover, its $\mu$ is the divergence of  (\ref{normH}).

\section{Results}

We say that a function $H\in C^1(A,\R)$, $A$ open subset of $\Omega$, is a {\it first integral} of (\ref{sysV}) on $A$ if $H$ is non-constant on any open subset of $A$, and, for every orbit $\gamma$, $I$ is constant on $\gamma \cap A$. By extension, we say that $H$ is a first integral of the vector field $V$.  

In next lemma we prove a relationship between normalizers and first integrals. We prove it for a period annulus, but it can be adapted to general systems admitting a first integral.

\begin{lemma}\label{prop}   Let $H$ be a first integral of (\ref{sysV}), and $A$ a period annulus of (\ref{sysV}). Assume $\nabla H$ not to vanish  on $A$.
Then a transversal vector field $W$ is a non-trivial normalizer of $V$ if and only if there exists a function $\xi\ne 0$ such that
$$
\partial_W H = \xi(H). 
$$
\end{lemma}
{\it Proof.} Let $W$ be a non-trivial normalizer of $V$. Let us choose arbitrarily a cycle $\gamma^*$ and a point $z^* \in\gamma^*$. Every point $z\in A$ can be written as $z=\phi_W(s,\phi_V(t,z^*))$. $W$ is a normalizer, hence  the parameter $s$ depends only on the cycle to which $z$ belongs. Hence the function that associates to a point $z\in A$ the value $s(z)$ of the parameter such that $z=\phi_W(s(z),\phi_V(t,z^*))$ is a first integral of (\ref{sysV}). By construction, one has
$$
\partial_W s = \nabla s W = 1.
$$
The above formula also implies that $\nabla s$ does not vanish on $A$. Hence there exists a scalar function $\chi$ such that $H(z) = \chi(s(z))$, with $\chi'(s) \neq 0$ because both $\nabla s$ and $\nabla H$ do not vanish. Then
$$
\partial_W H(z) = \chi' (s(z)) \partial_W s(z) =  \chi' (s(z))  = \chi'(\chi^{-1}(H(z))). 
$$
Then it is sufficient to set $\xi (H)= \chi' ( \chi^{-1}(H))$.  

Vice-versa, let us assume that there exists a scalar function $\xi$ such that $ \partial_W H = \xi(H) $. $W$ is transversal to $V$ because $V$ is orthogonal to $\nabla H$ and $ \partial_W H = \nabla H \cdot  W = \xi(H) \ne 0$. 
Moreover, since  $\nabla H$ does not vanish on $A$, $H$ does not has the same value on different cycles, so that every cycle in $A$ can be identified as $H^{-1}(l) \cap A$, for some $l\in \R$. This establishes a one-to-one correspondence between the cycles of $A$ and the values of $H$ on $A$. 
Let $z_1\neq z_2$ be distinct points of the same cycle $\gamma$. Then $H(z_1) = H(z_2)$. Since $\partial_W H = \xi(H)$, that is $H(\phi_W(s,z))$ depends only on the initial value of $H$ (in particular, it does not depend on the initial point $z$),  one has $H(\phi_W(s,z_1)) = H(\phi_W(s,z_2))$ for all $s$ for which both terms exist.  Hence the $W$-local flow takes arcs of orbits of (\ref{sysV}) into arcs of orbits of (\ref{sysV}), that is, $W$ is a normalizer of $V$. 
$\clubsuit$\bigskip

In order to compute $\mu$  for systems with a first integral, we proceed in two steps. First we work on hamiltonian systems, 
\begin{equation}\label{sysham}
x' = H_y,  \qquad y' = -H_x,
\end{equation}
then we extend our results to non-hamiltonian ones.  Let us  call  $V_H$ the vector field of (\ref{sysham}) and $W_H$ the vector field of (\ref{normH}). Such a system turns out to be a normalizer of (\ref{sysV}). In next lemma we prove that and compute its $\mu$. 
If $|\nabla H| \neq 0$ on $A$, every cycle is a level set of $H$, so that the period function is a function of the hamiltonian, that we denote by $T(H)$. 

\begin{lemma}\label{intnorm}   If $|\nabla H| \neq 0$ on $A$, then $V_H$  is a normalizer of  $W_H$  on $A$, with
\begin{equation}\label{muham1}
\mu_H =   \mbox{div }  W_H =
\frac{( H_{yy} - H_{xx})H_x^2 - 4H_{xy}H_xH_y  + (H_{xx}  - H_{yy})H_y^2}{|\nabla H|^4} ,
\end{equation}
and 
\begin{equation}\label{TpH}
T'(H) =  \int_0^{T(H)} \mu_H (\gamma(t))dt ,
\end{equation}
\end{lemma}
{\it Proof.} 
As proved in \cite{W}, the Lie brackets of a couple of vector fields $V$, $W$ satisfy the following formula
 $$ 
 [V,W] = \left( - \partial_W \ln (V \wedge W) +  \mbox{div } W
 \right) V +  \left(\partial_V \ln (V \wedge W) - \mbox{div } V \right) W  .
 $$
 One has $V \wedge W_H \equiv 1$ and $ \mbox{div } V \equiv 0$, so that 
 $$ 
 [V,W_H] =  \left( \mbox{div } W_H  \right) V := \mu_H V.
 $$
The final expression of $\mu_H $ is the outcome of standard differentiation operations. 

As for $T'(H)$, denoting by $\phi_{W_H}(s,z)$ a solution to (\ref{normH}), one has both   $\partial_{W_H} H =1$ and $\partial_{W_H} s =1$, so that $H$ and $s$ differ by a constant. Hence $T'(s) = T'(H)$.
$\clubsuit$\bigskip

Since $ \mu_H = \mbox{div } W_H$, both the shape of the $V$-orbits and the value of $|\nabla H |$ contribute to the sign of $ T'(H)$. In fact, in the limit case that $|\nabla H |$ be constant on a cycle $\gamma$, the only important feature is the curvature of $\gamma$ at its points, that determines the divergence of $W_H$. This is the case of hamiltonian systems whose hamiltonian function is a solution to the {\it eikonal equation}, $|\nabla H | \equiv 1$.  On the other hand, in the limit case that  $\gamma$ have constant curvature at all of its points, it is the value of $|\nabla H |$ to determine $ T'(H)$, as for the systems
$$
x' = y\rho(x^2+y^2) ,  \qquad y' = -x\rho(x^2+y^2) .
$$

\bigskip

In general, if $H$ is a first integral, then for every $\zeta\in C^0(\R,\R)$, $\zeta \neq 0$, the system \begin{equation}\label{normzeta}
x' = \frac{H_x}{|\nabla H|^2} \zeta(H),  \qquad y' = \frac{H_y}{|\nabla H|^2} \zeta(H),
\end{equation}
is a hamiltonian normalizer of (\ref{sysV}). In fact, if $\Phi = \int_0^s\frac{d\sigma}{\zeta(\sigma)}$, then the function
$J(z) = \Phi(H(z))$ is a hamiltonian function generating (\ref{normzeta}) as a  normalizer:
$$
\frac{J_x}{|\nabla J|^2}  = \frac{H_x \Phi'(H)}{|\nabla H|^2 \Phi'(H)^2} =  \frac{H_x}{|\nabla H|^2} \zeta(H).
$$

\bigskip

\begin{remark}Every non-trivial normalizer is a linear combination of a given non-trivial normalizer and $V$. Following proposition 2 in \cite{FGG}, one can easily compute the new normalizing function. In fact, every normalizer $N^*$ can be written as follows,
$$
N^* = \psi N + g V, 
$$
where $\psi$ is either a first integral of  (\ref{sysV}) or a constant, and $g$ is any function. In this case one can also write the new $\mu^*$:
$$
\mu^* = \psi \mu + \partial_V g.
$$
Then one has
$$
\int_0^T \mu^*(\phi_V(t,z))\ dt = \int_0^T \big( \psi \mu + \partial_V g \big)( \phi_V(t,z))\ dt = 
 \psi (\gamma) \int_0^T \mu( \phi_V(t,z))\ dt .
$$
On the other hand, it may occur that one of $\mu$, $\mu^*$ changes sign, while the other one does not, making easier to prove $T$'s monotonicity.
\end{remark}

\begin{remark} The formula (\ref{muham1}) can be also written as follows
$$
\mu = \frac{\Lambda_{H}}{|\nabla H|^4},
$$
where $\Lambda_{H}$ has been defined in \cite{S1}, example 5. 
\end{remark}

The normalizer  (\ref{normH}) is not necessarily that one providing the simplest possible $\mu$. For instance, if $H(x,y) = F(y) + G(x)$ our procedure gives  a $\mu$ of the form
\begin{equation}\label{Sabanorm}
\mu =  \frac{(F'' - G'')(F'^2 - G'^2)}{(F'^2 + G'^2)^2} ,
\end{equation} 
while a more convenient choice consists in taking the system
\begin{equation}\label{FGGnorm}
x' = \frac{G(x)}{G'(x)}, \qquad y' = \frac{F(y)}{F'(y)}
\end{equation}
on the set $F'(y)G'(x) \neq 0$ (see \cite{FGG}). For such a system, calling $W$ the normalizer, one has
$$
\partial_W H = G' \frac{G}{G'} + F' \frac{F}{F'} = H,
$$
so that the function $\xi$ of lemma \ref{prop} satisfies $\xi(r) =r$.
In this case the corresponding $\mu$ is given by
$$
\mu(x,y) =   \left(\frac{G(x)}{G'(x)} \right)' + \left( \frac{F(y)}{F'(y)}\right)' -1 = 
1  - \frac{G(x)G''(x)}{G'(x)^2} - \frac{F(y)F''(y)}{F'(y)^2}.
$$
A possibile advantage of our $\mu$, with respect to that one found in \cite{FGG}, is the possibility to compute the period's derivative also on period annuli where either $G'(x)$ or $F'(y)$ vanish. An example is given by  the system of example (7) in \cite{S1}
$$
x' = 2y, \qquad y' = - 4x(x^2-1),
$$
which has a period annulus encircling two centers and two homoclinic orbits. The denominators $F'(y)$ and $G'(x)$ vanish on every cycle encircling the two homoclinic orbits, while (\ref{Sabanorm}) does not.

\bigskip

We recall that a {\it reciprocal integrating factor (RIF)} of a differential system (\ref{sysV}), or of its vector field $V$, is a function $\kappa\in C^1(\Omega,\R)$, $\kappa > 0$,  such that its reciprocal $\frac{1}{\kappa}$ is an integrating factor of (\ref{sysV}). It is easy to prove that $\kappa$ is a RIF if and only if it satisfies $\kappa_x P + \kappa_y Q = \kappa (P_x + Q_y)$, that is $\partial_V \kappa = \kappa {\ \rm{div}} V$.

We report next lemma without proof. It shows that a normalizer of $V$ is also a normalizer of every reparametrization of $V$, and gives a formula for the new normalizing function. 

\begin{lemma}\label{normalizzatori} Assume $[V,W] = \mu V$, $\kappa > 0$ on $\Omega\setminus \{O\}$. Then $[\kappa V,U] =\overline \mu \left( \kappa V \right)$, with 
$$
\overline \mu = \mu -  \partial_W \ln \kappa. 
$$ 
\end{lemma}

If (\ref{sysV}) has a first integral $H$ with non-vanishing gradient, then we may write the system (\ref{sysV}) as a 
re-parametrized hamiltonian system,
$$
x' = P = H_y \kappa, \qquad y' = Q =  - H_x \kappa,
$$
where $\kappa$ satisfies
$$
\kappa = \frac{|V|}{|\nabla H |}.
$$

\begin{theorem}\label{rif}  Let  $\kappa$ a RIF of (\ref{sysV}). Then 
\begin{equation}\label{normkappa}
x' = -\kappa\ \frac{Q}{P^2+Q^2} = - \frac{Q}{|V| |\nabla H |}  , \qquad y' = \kappa\ \frac{P}{P^2+Q^2} = \frac{P}{|V| |\nabla H |}
\end{equation}
is a normalizer of $V$, and the function $\mu$ has the following form,
\begin{equation}\label{mukappa}
\mu =  \kappa \ \frac{- P^2(Q_x + P_y) + 2PQ (P_x - Q_y ) + Q^2 (Q_x + P_y) }{(P^2 + Q^2)^2} . 
\end{equation}
\end{theorem} 
{\it Proof.} The system (\ref{normkappa}) is just the system (\ref{normH}) written for $P = H_y \kappa$, $Q = -
H_x\kappa$.  The form of $\mu$ in (\ref{mukappa}) has been obtained by applying the formula
$$
\mu = \frac{[V,W]\cdot V}{|V|^2},
$$
which is a consequence of  $[V,W] = \mu V$.
$\clubsuit$\bigskip

Formula (\ref{mukappa}) has an evident relationship with the expression one gets for $\eta$ in  \cite{S1}, when choosing $V^\bot$ as the transversal vector field  (see \cite{S1}, corollary 4).  In fact, denoting such a function by $\eta^\bot$   one has
$$
\mu = \eta^\bot\ \frac{\kappa }{|V|^2} =  \frac{ \eta^\bot }{|V|  |\nabla H |} .
$$

\enddocument